\documentclass[11pt]{article}
\usepackage{amssymb}
\usepackage{graphicx}
\newtheorem {thm}{Theorem}[section]
\newtheorem {lem}[thm]{Lemma}

\def\ee{\epsilon}

\def\Cox{\hfill \Box}

\def\Z{{\mathbb Z}}

\def\br{\hfill \break}
\def\P{{\mathbb P}}
\def\E{{\mathbb E}}

\begin{document}

\begin{titlepage}
\begin{center}
{\large \bf Irreducible compositions \\ and the first return to the origin 
of a random walk} \\
\end{center}
\vspace{5ex}
\begin{flushright}
Edward A.\ Bender \footnote{University of California at San Diego, Department 
of Mathematics, 9500 Gilman Drive, Dept 0112, La Jolla, CA 92093-0112, 
ebender@ucsd.edu}\\
Gregory F.\ Lawler \footnote{Research supported in part by National Science 
Foundation grant \#  DMS 9971220}$^,$\footnote{Cornell University, 
Department of Mathematics, 310 Malott Hall, Ithaca, NY 14853-4201, 
lawler@math.cornell.edu}\\
Robin Pemantle \footnote{Research supported in part by
National Science Foundation grant \# DMS 0103635}$^,$\footnote{University
of Pennsylvania, Department of Mathematics, 209 S. 33rd Street, Philadelphia,
PA 19104-6395, pemantle@math.upenn.edu} \\
Herbert S.\ Wilf \footnote{University of Pennsylvania, Department of 
Mathematics, 209 S. 33rd Street, Philadelphia, PA 19104-6395, 
wilf@math.upenn.edu} 
\end{flushright}

\vfill

{\bf ABSTRACT:} \br 
Let $n = b_1 + \cdots + b_k = b_1' + \cdot + b_k'$ be a pair of
compositions of $n$ into $k$ positive parts.  We say this pair is 
{\em irreducible} if there is no positive $j < k$ for which $b_1 + \cdots b_j 
= b_1' + \cdots b_j'$.  The probability that a random pair of compositions
of $n$ is irreducible is shown to be asymptotic to $8/n$.  This problem
leads to a problem in probability theory.  Two players move along a game
board by rolling a die, and we ask when the two players will first coincide.
A natural extension is to show that the probability of a first return to 
the origin at time $n$ for any mean-zero variance $V$ random walk is
asymptotic to $\sqrt{V/(2 \pi)} n^{-3/2}$.  We prove this via two
methods, one analytic and one probabilistic.  
\vfill



\end{titlepage}

\setcounter{equation}{0}
\section{Introduction}

By a composition of $n$ into $k$ parts we mean an ordered representation 
of $n$ as a sum of $k$ positive integers. Let ${\cal C}$ and ${\cal C}'$
denote respectively the compositions $n = b_1 + \cdots + b_k$ and 
$n = b_1' + \cdots + b_k'$ of $n$ into $k$ parts. We'll say that 
${\cal C},{\cal C}'$ are an {\em irreducible pair} if for every 
$j=1, 2, \dots , k-1$ we have $b_1 + \cdots + b_j \neq b_1' + \cdots + b_j'$, 
while, of course, equality holds at $j=k$.  Note that we allow 
$b_1 + \ldots + b_i = b_1' + \cdots + b_j'$ for $i \neq j$.

Our starting point for this note is the following.  
Let $f(n)$ denote the number of irreducible ordered pairs 
of compositions of $n$ into the same number of parts.
\begin{thm} \label{thm:irreducible}
\begin{equation} \label{eq:comp GF}
\sum_{n \geq 1} f(n) z^n ~=~ \frac{z}{\sqrt{1 - 4z} + z} \, .
\end{equation} 
This gives a combinatorial interpretation to sequence 
\texttt{A081696} of Sloane's database.  Furthermore, if we let
$p_n$ denote the probability that a pair of compositions is
irreducible when chosen uniformly at random from among
all pairs of $n$ compositions into an equal number of parts, then
\begin{eqnarray}
f(n) & \sim & \frac{2}{\sqrt{\pi}} n^{-3/2} 4^n
   \, ; \label{eq:comp asym} \\[2ex] 
p_n & \sim & \frac{8}{n}
   \, . \label{eq:p}
\end{eqnarray}
\end{thm}

\noindent{\em Remark:} A somewhat similar problem about integer 
partitions was studied by Erd\H os et al~\cite{ENS}.

The asymptotics in~(\ref{eq:comp asym}) and~(\ref{eq:p})
are derived from the exact computation~(\ref{eq:comp GF}).
We compare this to another well known paradigm for analyzing
compositions, namely {\em poissonization}.  It is well known
that a uniform random composition of $n$ may be generated by
letting $\{ Y_j : j \geq 1 \}$ be independent random variables
whose distribution is geometric with mean~2, that is, $\P (Y_j = i)
= 2^{-i}$.  Let $W_k := \sum_{j=1}^k Y_j$ denote the partial sums.
If $T := \min \{ k : W_k \geq n \}$ is the first time the
partial sums exceed $n$, then $Y_1 + \cdots + Y_{T-1} + (n - Y_{T-1})$
is uniformly distributed over compositions of $n$.  The number of 
parts of the composition is $T$, which is asymptotically normal
with mean $n/2$ and standard deviation $\Theta (\sqrt{n})$.  
Conditioning on $T=k$ gives the uniform distribution on compositions 
of $n$ into $k$ parts.  

 From this viewpoint, a pair of compositions with the same number of
parts is just a pair of independent random walk sequences 
$\{ Y_j \}$ and $\{ Y_j' \}$, conditioned to have the same 
stopping time $T = T'$.  Irreducibility of the pair corresponds
to $W_j \neq W_j'$ for all $1 \leq j \leq T-1$.  Let 
$S_k := W_k - W_k' = \sum_{j=1}^k (Y_j - Y_j') := \sum_{j=1}^k X_j$ 
be the partial sums of the difference sequence $\{ X_j \}$.  
Then irreducibility corresponds to $\tau = n$, where $\tau$ is the 
first return time, that is, $\tau = \min \{ k \geq 1 : S_k = 0 \}$.  

A rigorous proof of~(\ref{eq:p}) via analysis of the return
time of $\{ S_n \}$ to the origin would require, among other things, 
showing that conditioning on $T = T'$ does not significantly
affect the distribution of the return time.  This would be far
messier than the compact proof of Theorem~\ref{thm:irreducible}
below.  Nevertheless, the poissonization
paradigm raises the question of the distribution of the return time 
of $\{ S_n \}$ to the origin.  The same question arises in a game
similar to Parcheesi with only one token per player.  Here, the 
two players each roll a die and (simultaneously) advance their 
token the number of positions shown on the die.  When the tokens 
collide, they must both go back to start.  The chance of 
the first collision occurring at time $n$ is just
$a_n := \P (\tau = n)$.  Our main result is the following
asymptotic for $a_n$:
\begin{thm} \label{th:finite variance}
Let $\{ X_j : j \geq 1 \}$ be independent with mean zero, 
finite variance, $V$, and no periodicity (that is, the GCD
of times $n$ at which it is possible to have $S_n = 1$ is~1).
Let $\{ S_n \}$, $\tau$ and $a_n$ be as 
above.  Then the probability $a_n$ of the first return to~0 
occurring at time $n$ is asymptotically given by
$$a_n ~\sim~ \sqrt{\frac{V}{2 \pi}} n^{-3/2} \, .$$
\end{thm}

Surprisingly, given the wealth of knowledge about random walks,
we were unable to find this theorem in the literature.
The formula is not surprising, and is what one obtains 
in a thumbnail calculation by ``differentiating'' with respect
to $n$ estimates such as~(\ref{eq:Q}) below, which is an estimate
for $Q_n := \P(S_j \neq 0, \forall 1 \leq j \leq n)$. 
There are special cases, such as the simple random walk where 
$S_n = \pm 1$ according to a fair coin-flip, in which $a_n$ 
is easy to compute exactly.  Asymptotics in the general case
are well known for many quantities such as $Q_n$ and 
$a_n' := \P (S_n = 0)$, but we could find no text that 
included asymptotics for $a_n$ and indeed these seem tricky 
to obtain by probabilistic methods; a probabilistic proof of 
Theorem~\ref{th:finite variance} is the subject of the 
last section of this note.

In the remainder of this section, we prove Theorem~\ref{thm:irreducible}.
In the subsequent section we prove Theorem~\ref{th:finite variance}
by analytic means.  In the final section, we give a probabilistic 
proof of Theorem~\ref{th:finite variance}.  

\noindent{\sc Proof of Theorem}~\ref{thm:irreducible}:
Let $f(n,k)$ be the number of irreducible ordered pairs of
compositions of $n$ into $k$ parts.  We will show that
\begin{equation}
\label{eq:twovar}
\sum_{n , k\geq 1}f(n,k) x^n y^k ~=~ 
   \frac{xy \left ( \sqrt{1+x^2(1-y)^2-2x(1+y)} - xy \right ) }
   {1-2x(1+y)+x^2(1-2y)} \, ,
\end{equation}
from which (\ref{eq:comp GF}) follows by setting $y = 1$.

To show (\ref{eq:twovar}), by considering the number of ordered pairs of
compositions of $n$ into $k$ parts such that the partial sums of the parts
agree with each other at indices $k_1 , k_1 + k_2 , \dots , k_1 + \dots 
+ k_r$, we see that
$$\sum_{r \geq 1} \sum_{{m_1+\dots +m_r=n}\atop{k_1+\dots +k_r=k}}
   f(m_1 , k_1) f(m_2 , k_2) \dots f(m_r , k_r) ~=~ {n-1\choose k-1}^2 \, ,$$
the right side being the total number of pairs of compositions of $n$ into
$k$ parts. Hence if $F(x,y) = \sum_{n , k \geq 1} f(n,k)x^ny^k$, we have
\begin{eqnarray*}
\frac{F}{1-F} = F+F^2+F^3+\dots & = & \sum_{n , k \geq 1}
   {n-1\choose k-1}^2 x^n y^k \\[1ex]
& = & xy \sum_{n , k \geq 0}{n\choose k}^2 x^n y^k \\[1ex]
& = & xy \sum_{n \geq 0} x^n (1-y)^n P_n \left ( \frac{1+y}{1-y} \right )
   \\[1ex]
& = & \frac{xy}{\sqrt{1-2x(1+y)+x^2(1-y)^2}} \, , 
\end{eqnarray*}
in which the $P_n$'s are the Legendre polynomials. The claimed result
(\ref{eq:twovar}) now follows by solving for $F$.

The estimate~(\ref{eq:comp asym}) follows from~(\ref{eq:comp GF})
via standard Tauberian theorems.  The result of Flajolet and Odlyzko,
for instance (Theorem~\ref{th:FO} quoted below) suffices,
although~(\ref{eq:comp asym}) may also be obtained by the method
of Darboux which requires a smaller region of analyticity.
Since there are
$$ \sum_k {n-1\choose k-1}^2 ~=~ {2n-2\choose n-1} \sim
   \frac{4^{n-1}}{\sqrt{n\pi}} $$
ordered pairs of compositions of $n$ with the same number of parts, it
follows that the probability that a random pair of compositions of $n$ with
the same number of parts is irreducible is $\sim 8/n$.   $\Cox$

\section{Analytic proof}

Let $H(z) := \sum_{n \geq 1} a_n z^n$ be the generating function for
the probabilities $a_n$ of first return at time $n$.  Let $G(z) := 
\sum_{n \geq 0} a_n' z^n$, where $a_n' = \P (S_n = 0)$ is
the probability of a return to the origin at time $n$ but not
necessarily the first return (set $a_0' = 1$ and $a_0 = 0$).   
Then $G$ and $H$ are analytic on the open unit disk and 
$G = 1 / (1 - H)$.  We will use this to obtain $H$ from $G$,
while $G$ in turn is obtained from the two-variable generating function
$$F(z,w) ~:=~ \sum_{n \geq 0} \sum_{j \in \Z} \P (S_n = j) z^n w^j \, .$$
Finally, we may write $F = 1 / (1 - z g(w))$ where 
$$g(w) ~:=~ \sum_{n \in \Z} b_n w^n$$
is the generating function for $X_1$.  

The following estimates are elementary.  From the local central
limit theorem~\cite[Theorem~(II.5.2)]{Dur}, as $n \to \infty$,
\begin{equation} \label{eq:CLT}
a_n' ~\sim~ \frac{1}{\sqrt{2 \pi V}} n^{-1/2} \, .
\end{equation}
Consequently,
\begin{equation} \label{eq:G}
G(z) ~\sim~ \frac{1}{\sqrt{2V}} (1-z)^{-1/2}
\end{equation}
as $z \uparrow 1$.  To see this, let $\ee$ denote $1 - z$ and compute
\begin{eqnarray*}
G(z) & = & \sum_{n \geq 0} \frac{1}{\sqrt{2 \pi V}} n^{-1/2} 
      e^{-n \ee (1 + o(1))}  \\
& = & \frac{1}{\sqrt{2 \pi V}} \ee^{1/2} \sum_{n \geq 0} (n \ee)^{-1/2} 
      e^{- n \ee (1 + o(1))} \\
& \sim & \frac{1}{\sqrt{2 \pi V}} \ee^{-1/2} \int_0^\infty x^{-1/2} e^{-x}
   ~=~ \frac{1}{\sqrt{2 V}} (1-z)^{-1/2}
\end{eqnarray*}
using dominated convergence at the first approximation.

Finally, for $H = 1 - 1/G$, we have the estimate
\begin{equation} \label{eq:H}
1 - H (z) ~\sim~ \sqrt{2V} (1-z)^{1/2}
\end{equation}
as $z \uparrow 1$. 
The proof of Theorem~\ref{th:finite variance} rests on these estimates
and on the following Tauberian theorem of~\cite{FO}:
\begin{thm}[Flajolet-Odlyzko (1990)]  \label{th:FO}
Say that a region $R$ is a {\em Camembert region} \footnote{Named, 
by French mathematicians, for its shape.} if it is of the form
$R_\ee := \{ |z| < 1 + \ee ~\mbox{ and }~ |\arg(z-1)| > \pi/2 - \ee \}$. 
If a function $H$ is analytic in a Camembert region and
$H(z) \sim C (1-z)^{- \alpha}$ near $z = 1$, then its coefficients 
$a_n$ satisfy
$$ a_n ~\sim~ \frac{C}{\gamma (\alpha)} n^{\alpha - 1} \, .$$
$\Cox$
\end{thm}

\noindent{\sc Proof of Theorem}~\ref{th:finite variance}:
Let $C$ denote the unit circle.  For fixed $w \in C$, the 
function $F(z,w)$ is analytic as $z$ varies over the open unit disk; 
this follows from absolute convergence of the power series.  
It also follows that $F(z,w)$ is continuous in $(z,w)$ on
the product $\Omega := D \times C$ of the open unit disk 
with the unit circle.  For fixed $w \in C$, the Cauchy integral 
formula gives
$$z^n \P (S_n = 0) ~=~ z_n \int \sum_j \P (S_n = j) w^n \frac{dw}{w} \, .$$
We may sum this over $n$ and exchange the sum and integral
as long as $|z| < 1$, leading to
\begin{equation} \label{eq:cauchy}
G(z) ~=~ \frac{1}{2 \pi i} \int_\gamma \frac{F(z,w)}{w} \, dw
\end{equation}
where $\gamma$ goes around the unit circle, counterclockwise\footnote{ 
This integral formula is used in~\cite{HK} to derive a result
(attributed to~\cite{Fur} by~\cite{Sta}) implying in this case
that $G$ is algebraic whenever $g$ is rational.  In fact, in the 
case where $X_1$ has finite support, one may use this implication
at the next step to avoid having to examine the power series expansion 
of $g$.}.

Suppose we can show $G$ to be analytic in a Camembert region.  
It follows that $H = 1 - 1/G$ is meromorphic in a Camembert region,
and since a function whose coefficients go to zero may have no
poles in the closed unit disk, it follows that $H$ is analytic
in a Camembert region.  The conclusion of the theorem will then 
follow from~(\ref{eq:H}) and Flajolet-Odlyzko.

\begin{quote} 
{\bf Claim:} There is a Camembert region $R$ such that $z \neq 1 / g(w)$ 
for any $z \in R$ and $w$ on the unit circle, $C$.  Consequently,
$F(z,w)$ has an extension to $R \times C$ that is analytic in $z$
and continuous in $(z,w)$.

\noindent{\sc Proof:} The facts that $X_1$ is a probability distribution,
has mean zero, and has variance $V$ translate into three facts
about $g$, namely, $g(1) = 1, g'(1) = 0, g''(1) = V$.  Immediately,
we then have 
$$\frac{1}{g(e^{i \theta})} ~=~ 1 + \frac{V}{2} \theta^2 + o(\theta^2) \, .$$
Hence $\arg (1/g(w)) \to 0$ as $w \to 1$ in $C$ and there is an
$\delta > 0$ such that for $\arg (z-1) > \delta$ and $|\theta| < \delta$, 
$z \neq 1 / g(e^{i \theta})$.  For $1 \neq w \in C$, aperiodicity of 
$X_1$ implies $|g(w)| < 1$.  Let $\ee$ be the minimum of $\delta$ and
the values $|g(e^{i \theta})|^{-1} - 1$ on $|\theta| \in [-\pi , \pi] 
\setminus (-\delta , \delta)$.  Then $z \neq 1 / g(w)$ on the 
Camembert region $R(\ee)$.   $\Cox$
\end{quote} 

Finishing the proof of Theorem~\ref{th:finite variance}, we
let $R$ be as in the conclusion of the lemma and observe that 
for any closed loop $\beta$ in $R$, we may exchange the order of
integration in the representation of $G$ in~(\ref{eq:cauchy}) to get 
\begin{eqnarray*}
\int_\beta G(z) \, dz & = & \int_\beta \int_\gamma 
   \frac{1}{2 \pi i} \frac{F(z,w)}{w} \, dw \, dz \\[1ex]
& = & \frac{1}{2 \pi i} \int_\gamma \, \frac{dw}{w} \, 
   \int_\beta F(z,w) \, dz \\[1ex]
& = & 0 \, .
\end{eqnarray*}
By Morera's theorem, $G$ is analytic in $R$, completing the proof
of Theorem~\ref{th:finite variance}.
$\Cox$  

\section{Probabilistic proof}

Let $p_n (x,y) := \P_x (S_n = y)$ and $q_n (x,y) := \P_x (S_n = y , 
S_j \neq 0 \; \forall 1 \leq j \leq n-1)$ denote probabilities for 
$\{ S_n = y \}$ respectively with or without killing at the origin.  
These quantities are symmetric in the two arguments.  Previously
defined quantities are related to these by $a_n' = p_n (0,0)$ 
and $a_n = q_n (0,0)$.  We let $Q_n := \P_0 ( S_j \neq 0 , \, \forall 
1 \leq j \leq n) = \sum_{k=n+1}^\infty a_n$ be the tail sums of $\{ a_n \}$.

The derivations of~(\ref{eq:CLT}),~(\ref{eq:G}) and~(\ref{eq:H}) in
the probability literature are via the generating function analysis
in the previous section.  At this point, the methods part ways.  
The probabilistic analysis derives $a_n$ from its tail sums, $Q_n$.
The reasonably well known estimate on $Q_n$ is
\begin{equation} \label{eq:Q}
Q_n ~\sim~ \sqrt{\frac{2V}{\pi}} n^{-1/2} \, .
\end{equation}
This is proved analytically, not via extending a two-variable 
generating function to a Camembert region, but just from~(\ref{eq:H}). 
The key here is that the sequence $\{ Q_n \}$ is monotone.  
According to a Tauberian theorem which may be found 
in~\cite[Theorem~XIII.5]{Fel}, the extra regularity, together 
with the behavior of its generating function $H / (1-z)$ for real 
$z \uparrow 1$, implies~(\ref{eq:Q}).

Since $a_{n+1} = Q_n - Q_{n+1}$, the conclusion of the theorem now
follows if we can establish regularity of $a_n$ to the degree that
\begin{equation} \label{eq:regular}
a_{n+1} - a_n ~=~ O(n^{-5/2}) \, .
\end{equation}
Note that we have now converted the task from one of finding the
correct leading term into one of finding an upper bound to within
a constant factor, which is a problem well suited to probabilistic 
analysis.  To complete the regularity argument we need a couple of
estimates on how rapidly $p_n (0,x) := \P_0 (S_n = x)$ can change with $n$.
We will prove these at the end.  
\begin{lem} \label{lem:bounds}
Under the assumptions of aperiodicity, zero mean and finite variance,
there is a constant $C$ such that
\begin{eqnarray}
\left | p_n (0,x) - p_{n+1} (0,x) \right | & \leq & \frac{C}{n^{3/2}} 
   \label{eq:small p} \\[1ex]
\left | p_n (0,x) - p_{n+1} (0,x) \right | & \leq & \frac{C}{(1+x^2) n^{1/2}} 
   \label{eq:large p}  \, .
\end{eqnarray}
Upper bounds for $q_n (0,x)$ are given by
\begin{eqnarray} \label{eq:q}
q_n (0,x) & \leq & \frac{C (|x| + 1)^{1/2}}{n^{3/2}} \label{eq:small q} 
   \\[1ex]
q_n (0,x) & \leq & \frac{C}{n} \label{eq:large q}\, .
\end{eqnarray}
\end{lem}

We now prove~(\ref{eq:regular}) for $n = 3m$, the cases of $3m+1$ and
$3m+2$ being identical.  Breaking down according to location at
times $m$ and $2m$ we get
$$a_n ~=~ q_{3m} (0,0) ~=~ \sum_{x \neq 0} \sum_{y \neq 0} 
   q_m (0,x) q_m (x,y) q_m (y,0) \, .$$
Write $q_m (x,y) = p_m (x,y) - \P_x (S_m = y , \exists j \in [1,m-1] 
\, : \, S_j = 0)$.  
Substituting this in the above equation gives
\begin{eqnarray*}
a_n & = & \sum_{x \neq 0} \sum_{y \neq 0} q_m (0,x) p_m (x,y) q_m (y,0) 
   - \sum_{k=1}^{m-1} \sum_{y \neq 0} a_{m+k} p_{m-k} (0,y) q_m (y,0)
   \\[1ex]
a_{n+1} & = & \sum_{x \neq 0} \sum_{y \neq 0} q_m (0,x) p_{m+1} (x,y) q_m (y,0) 
   - \sum_{k=1}^m \sum_{y \neq 0} a_{m+k} p_{m+1-k} (0,y) q_m (y,0) 
\end{eqnarray*}
and hence
\begin{eqnarray*}
|a_n - a_{n+1}| & \leq & \sum_{x \neq 0} \sum_{y \neq 0} 
   q_m (0,x) |p_m (x,y) - p_{m+1} (x,y)| q_m (y,0) \\
&& + \sum_{y \neq 0} a_{2m} p_1 (0,y) q_m(0,y) \\
&& + \sum_{k=1}^{m-1} \sum_{y \neq 0} a_{m+k} |p_{m+1-k} (0,y) - 
   p_{m-k} (0,y)| q_m (y,0) \, .
\end{eqnarray*}

We must bound each of the three terms by $O(m^{-5/2})$.  The second term 
is $a_{2m} a_m = O(m^{-3})$ by~(\ref{eq:small q}).  Using~(\ref{eq:small p})
we see that the first term is 
$$O(m^{-3/2}) \sum_{x \neq 0} q_m (0,x) \sum_{y \neq 0} q_m (y,0) 
   ~=~ O(m^{-3/2}) Q_m^2 ~=~ O(m^{-5/2}) \, .$$
The third term requires a little more care.  We will show that
\begin{equation} \label{eq:summand}
\sum_{y \neq 0} |p_{k+1} (0,y) - p_k (0,y)| q_m (0,y) ~\leq~ c k^{-1/2} 
   \left [ 1 + \log \left ( \frac{m}{k} \right ) \right ] m^{-3/2} \, .
\end{equation}

To show this, split into three ranges of values for $y$, namely
$|y| \leq \sqrt{k}$, $\sqrt{k} < |y| <  \sqrt{m}$ and $|y| \geq \sqrt{m}$.
In the first range we use~(\ref{eq:small p}) and~(\ref{eq:small q}) 
with $|y|$ bounded above by $k^{1/2}$
to see that the summand is $O(k^{-3/2}) O(k^{1/2} m^{-3/2})$.  There
are $k^{1/2}$ summands, so the total sum is $O(k^{-1/2} m^{-3/2})$. 

In the middle range, we use~(\ref{eq:large p}) and~(\ref{eq:small q})
to see that the summand is bounded by a constant multiple of
$k^{-1/2} |y|^{-2} |y| m^{-3/2}$.  Summing over $y$ introduces
the factor of $(1/2) \log (m/k)$.  For the third sum, use~(\ref{eq:large p})
and~(\ref{eq:large q}) to see that the summand is 
$O(k^{-1/2} |y|^{-2} m^{-1})$, so that summing over $y \geq \sqrt{m}$
gives $O(k^{-1/2} m^{-3/2})$.  This proves~(\ref{eq:summand}).  

Finally, summing~(\ref{eq:summand}) over $k < m$ gives $O(m^{-5/2})$
which establishes~(\ref{eq:regular}), finishing the proof of
Theorem~\ref{th:finite variance}.   $\Cox$

\noindent{\sc Proof of Lemma}~\ref{lem:bounds}:  
The simplest of the inequalities is~(\ref{eq:large q}), so we
handle it first.  Let $n = 3m$.  The bound~(\ref{eq:large q}) 
follows immediately from
$$q_n (0,x) ~\leq~ \sum_y q_m (0,y) p_{2m} (y,x) 
   ~\leq~ a_m \sup_{y,x} p_{2m} (y,x) ~=~ O(m^{-1}) \, .$$

To prove~(\ref{eq:small q}), we decompose according to the
position at time $m$ and at time $2m$, so that
\begin{eqnarray*}
q_n (0,x) & = & \sum_{y,z \neq 0} q_m (0,y) q_m (y,z) q_m (z,x) \\
& \leq & \left ( Q_m \sup_{y,z} q_m (y,z) \right ) \sum_{z \neq 0} 
   q_m (z,x) \\
& = & O(m^{-1}) Q_m (x) \, ,
\end{eqnarray*}
where $Q_m (x):= \P_x (S_j \neq 0 \; \forall 1 \leq j \leq m$ and we 
have used $q_m (z,x) = q_m (x,z)$ to infer $\sum_z q_m (z,x) = Q_m (x)$.  
Thus it suffices to show that 
\begin{equation} \label{eq:Q_m}
Q_m (x) ~=~ O((1 + |x|) m^{-1/2}) \, .
\end{equation}

Observe that there is a constant $c$
depending on the distribution of $X_1$ but not on $y$ such
that the probability, call it $\rho (y)$, of hitting $y$ in at
most $y^2$ steps starting from the origin is at least $c$ 
(use the local central limit theorem
to bound the expected number of visits to $y$ within the first
$t^2$ steps by from below by $c_1 (1+|y|)$ and use the Green's function
to bound the expected number of visits to $y$ given at least
one visit from above by $c_2 (1+|y|)$).  By a last exit decomposition,
we then have
$$c ~\leq~ \rho (y) ~\leq~ \sum_{j \leq y^2} p_j (0,0) \tau (y)$$
where $\tau (y)$ is the probability starting at the origin of
hitting $y$ before returning to the origin.  Using $p_j (0,0) = 
\Theta (j^{-1/2})$ and solving for $\tau (y)$ gives 
$$\tau (y) ~=~ \Omega \left ( \frac{1}{1 + |y|} \right ) \, .$$
But by~(\ref{eq:Q}), decomposing according to the time $y$ is first hit,
$$\sqrt{\frac{2V}{\pi}} n^{-1/2} ~\sim~ Q_n \geq \tau (y) Q_n (y)$$
and solving for $Q_n(y)$ proves~(\ref{eq:Q_m}).  

The bounds on $\delta_n (x) := |p_n (0,x) - p_{n+1} (0,x)|$ are 
classical (though not all that well known) and are obtained by
the same means as the local central limit theorem.  Let 
$$\phi (\theta) ~=~ g(i \theta) = \E e^{i \theta X_1}$$
be the characteristic function of $X_1$, so that as we have seen,
mean zero, finite variance and aperiodicity imply that
\begin{equation} \label{eq:phi 1}
\phi (\theta) ~=~ 1 - \frac{V \theta^2}{2} + o(\theta^2)
\end{equation}
near~1, while 
\begin{equation} \label{eq:phi 2}
1 - \phi (\theta) ~\leq~ c \theta^2
\end{equation}
for all $|\theta| \leq \pi$.  The inversion formula gives
\begin{equation} \label{eq:p1}
p_n (0,x) ~=~ \frac{1}{2 \pi} \int_{-\pi}^\pi 
   \phi (\theta)^n e^{- i \theta x} \, d\theta \, .
\end{equation}
We obtain from~(\ref{eq:p1})
$$\delta_n (x) ~\leq~ \int_{-\pi}^{\pi} |\phi (\theta)|^n 
   |\phi (\theta) - 1| d \theta \, .$$
Given~(\ref{eq:phi 1}) and~(\ref{eq:phi 2}), we see this is a 
saddle point integral with main contribution near $\theta = 0$. 
In particular, since $\phi (\theta) \leq 1 - c \theta^2$, we know
that $|\phi (\theta)|^n \leq c \exp (b n^{3/4})$ for $n^{-1/8}
\leq |\theta| \leq \pi$ and we may restrict the integrals to
a suitable range such as $|\theta| < n^{-1/8}$.  
Let $y = \theta \sqrt{n}$.  From~(\ref{eq:phi 1}),
\begin{eqnarray*}
\left | \phi \left ( \frac{y}{\sqrt{n}}\right ) \right |^n 
   & \leq & c e^{-a y^2} \, , \\
\left | \phi \left ( \frac{y}{\sqrt{n}} \right ) - 1 \right | 
   & \leq & c \frac{y^2}{n} \, , 
\end{eqnarray*}
whence $\delta_n (x) = O(n^{-3/2})$.  

For~(\ref{eq:large p}) we integrate~(\ref{eq:p1}) twice by parts to get
\begin{equation} \label{eq:p2}
p_n (0,x) ~=~ - \frac{x^2}{2 \pi} \int_{-\pi}^{\pi} 
   \left [ n(n-1) \phi'(\theta)^2 + n \phi'' (\theta) \right ] 
   \phi (\theta)^{n-2} e^{- i x \theta} \, d\theta \, .
\end{equation}
The same truncation and change of variables, together with the
estimate 
$$ \left | \phi' \left ( \frac{y}{\sqrt{n}} \right ) \right | 
   ~\leq~ c \frac{1 + |y|}{\sqrt{n}}$$ 
give $\delta_n (x) = O((1+x^2) n^{-1/2})$ and completes the proof 
of the lemma.   $\Cox$

\noindent{\bf Acknowledgement:} The problems we address
were suggested by a question posed by Dr.\ Amy Myers.

\end{document}